\documentclass{amsart}

\usepackage{graphicx}
\usepackage{psfrag}
\begin{document}

\newcommand{\pic}[2]{\includegraphics[scale=0.#1]{#2.eps}}
\newcommand{\D}{{\mathfrak D}}
\newcommand{\N}{{\mathbb N}}
\newcommand{\Z}{{\mathbb Z}}
\newcommand{\C}{{\mathbb C}}
\newcommand{\Schub}{{\mathfrak{S}}}
\newcommand{\Groth}{{\mathfrak{G}}}
\newcommand{\SL}{\operatorname{SL}}
\newcommand{\Hecke}{{\mathcal H}}
\newcommand{\wt}{\widetilde}
\newcommand{\cross}{\includegraphics[scale=0.1]{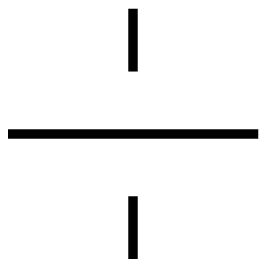}}
\newcommand{\avoid}{\includegraphics[scale=0.1]{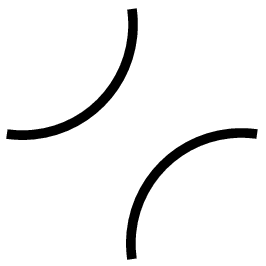}}
\newcommand{\comment}[1]{}

\newtheorem*{thm}{Theorem}
\newtheorem{cor}{Corollary}

\psfrag{y1}{$b_1$}
\psfrag{y2}{$b_2$}
\psfrag{y3}{$b_3$}
\psfrag{y4}{$b_4$}
\psfrag{y5}{$b_5$}
\psfrag{y6}{$b_6$}
\psfrag{y7}{$b_7$}
\psfrag{y8}{$b_8$}

\title{Specializations of Grothendieck polynomials}
\author{Anders S.~Buch}
\address{Matematisk Institut, Aarhus Universitet, Ny Munkegade, 8000
  {\AA}rhus C, Denmark}
\email{abuch@imf.au.dk}
\author{Rich\'ard Rim\'anyi}
\address{Department of Mathematics, The University of North Carolina
  at Chapel Hill,\linebreak CB \#3250, Phillips Hall, Chapel Hill, NC
  27599, USA}
\email{rimanyi@email.unc.edu}
\date{\today}


\maketitle

\section{Introduction}

Let $v, w \in S_n$ be permutations and let $\Schub_w(x;y)$ and
$\Groth_w(a;b)$ denote the double Schubert and Grothendieck
polynomials of Lascoux and Sch\"utzenberger
\cite{lascoux.schutzenberger:structure}.  The goal of this note is to
prove a formula for the specializations of these polynomials to
different rearrangements of the same set of variables.  For example
\[ \Groth_w(b_v;b) = \Groth_w(b_{v(1)},\dots,b_{v(n)} \,;\,
   b_1,\dots,b_n) \,.
\]

The double Schubert polynomial $\Schub_w(x;y)$ represents the class of
the Schubert variety for $w$ in the torus-equivariant cohomology of
$\SL_n(\C)/B$.  The specialization $\Schub_w(y_v;y)$ gives the
restriction of this class to the fixed point corresponding to $v$
\cite[Thm.~2.3]{goldin:cohomology}.  Equivalently, $\Schub_w(x;y)$
represents the class of an orbit in the $B \times B$-equivariant
cohomology of $\C^{n \times n}$ \cite{feher.rimanyi:schur,
  knutson.miller:grobner}, and $\Schub_w(y_v;y)$ is the restriction of
this class to another orbit, i.e. an `incidence class' in the sense of
\cite{rimanyi:thom}.  Specialized Grothendieck polynomials
$\Groth_w(b_v;b)$ have similar interpretations in equivariant
$K$-theory.

The formula proved in this paper generalizes the usual formulas for
Schubert and Grothendieck polynomials in terms of RC-graphs
\cite{fomin.kirillov:yang-baxter, fomin.kirillov:grothendieck,
  bergeron.billey:rc-graphs, knutson.miller:subword}, and it
furthermore gives immediate proofs of several important properties of
these polynomials.  This includes Goldin's characterization of the
Bruhat order \cite{goldin:cohomology}, the existence and supersymmetry
of stable Schubert and Grothendieck polynomials
\cite{fomin.kirillov:yang-baxter}, as well as the statements about
Schubert and Grothendieck polynomials needed in \cite{warsawa} and
\cite{buch:alternating}.  The proof of our formula relies on Fomin and
Kirillov's construction of Grothendieck polynomials
\cite{fomin.kirillov:grothendieck}.

We thank W.~Fulton and S.~Billey for helpful comments and
references, and L.~Feh\'er for inspiring collaboration on related
papers.

\section{The main theorem}

Consider the diagram $\D_v$ consisting of hooks of lines going due
north and due west from the points $(v(j),j)$, and with each such hook
labeled by $b_{v(j)}$.  For example, when $v = 264135$ we get:
\[ \D_v \ \ = \ \ \ \raisebox{-40pt}{\pic{70}{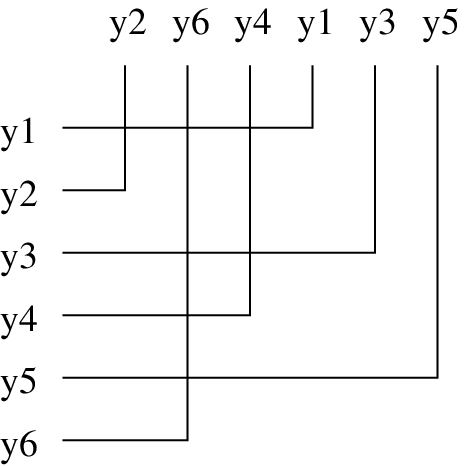}} \]
We let $C(\D_v)$ denote the crossing positions of this diagram, i.e.\
the points $(i,j)$ such that $v(j) > i$ and $v^{-1}(i) > j$.  Notice
that $C(\D_v) = D(v^{-1})$ with the notation of
\cite[p.~8]{macdonald:notes}.  For $(i,j) \in C(\D_v)$ we let
$\nu(i,j)$ be one plus the number of hooks going north-west of $(i,j)$
in the diagram $\D_v$, i.e.
\[ \nu(i,j) = j + \# \{ k > j : v(k) < i \} \,. \]

We need the {\em degenerate Hecke algebra}, which is the free
$\Z$-algebra $\Hecke$ generated by symbols $s_1, s_2, \dots$, modulo
the relations (i) $s_i s_j = s_j s_i$ if $|i-j| \geq 2$, (ii) $s_i
s_{i+1} s_i = s_{i+1} s_i s_{i+1}$, and (iii) $s_i^2 = -s_i$.  This
algebra has a basis of permutations.

For a subset $D \subset C(\D_v)$, consider the product in $\Hecke$ of
the simple reflections $s_{\nu(i,j)}$ for $(i,j) \in D$, in south-west
to north-east order, {i.e.\ }$s_{\nu(i,j)}$ must come before
$s_{\nu(i',j')}$ if $i \geq i'$ and $j \leq j'$.  This product is
equal to plus or minus a single permutation $w(D)$.  We say that $D$
is a {\em Fomin-Kirillov graph\/} (or {\em FK-graph\/}) for this
permutation w.r.t.\ the diagram $\D_v$, and that $D$ is {\em
  reduced\/} if $|D|$ equals the length of $w(D)$.

An FK-graph $D$ can be pictured by replacing the crossing positions of
$\D_v$ which belong to $D$ with the symbol ``\,$\cross$\,'', while the
remaining crossing positions are replaced with the symbol
``\,$\avoid$\,''.  If $D$ is reduced then the string entering the
resulting diagram at column $j$ at the top will exit at row $w(D)(j)$
at the left hand side.


Our main result is the following theorem, which is proved
combinatorially in the next section.  It is natural to ask for a
geometric proof as well.

\begin{thm}
For permutations $v,w \in S_n$ and variables $b_1,\dots,b_n$ we have
\[ \Groth_w(b_v; b) \ = \ \sum_D \, (-1)^{|D|-\ell(w)} \prod_{(i,j)\in D} 
   \left( 1 - \frac{b_i}{b_{v(j)}} \right)
\]
where the sum is over all FK-graphs $D$ for $w$ {w.r.t.\ }$\D_v$.
\end{thm}

\begin{cor} \label{cor:schubspec}
For permutations $v,w \in S_n$ and variables $y_1,\dots,y_n$ we have
\[ \Schub_w(y_v; y) \ = \ \sum_D \prod_{(i,j) \in D} (y_{v(j)} - y_i) \]
where the sum is over all reduced FK-graphs $D$ for $w$ {w.r.t.\ 
  }$\D_v$.
\end{cor}

\begin{cor}
  The usual formulas for double Schubert and Grothendieck polynomials
  in terms of RC-graphs are true (see
  \cite{fomin.kirillov:yang-baxter, fomin.kirillov:grothendieck,
    bergeron.billey:rc-graphs, knutson.miller:subword}).
\end{cor}
\begin{proof}
  Apply the theorem to $\Groth_w(a_1,\dots,a_n,b_1,\dots,b_n \,;\,
  b_1,\dots,b_n,a_1,\dots,a_n)$.
\end{proof}

The next corollary recovers the characterization of the Bruhat order
proved in \cite[Thm.~2.4]{goldin:cohomology}.

\begin{cor} \label{cor:bruhat}
Let $v,w \in S_n$.  The following are equivalent:
\begin{enumerate}
\item $w \leq v$ in the Bruhat order.
\item $\Schub_w(y_v; y) \neq 0$.
\item $\Groth_w(b_v; b) \neq 0$.
\end{enumerate}
\end{cor}
\begin{proof}
  The product defining $w(C(\D_v))$ is a reduced expression for $v$.
  There exists a reduced FK-graph $D \subset C(\D_v)$ for $w$ if and
  only if $w$ equals a reduced subexpression of this product.  The
  later is equivalent to $w \leq v$.
  This shows that each of (2) and (3) imply (1) (these implications
  are clear from geometry, too.)  It is also clear that (2) implies
  (3).  To see that (1) implies (2), notice that if $(i,j) \in
  C(\D_v)$ then $i < v(j)$.  Therefore each reduced FK-graph $D$ for
  $w$ in Corollary \ref{cor:schubspec} contributes a positive
  polynomial in the variables $z_i = y_{i+1} - y_i$.
\end{proof}

The following corollary implies that stable double Schubert and
Grothendieck polynomials exist and are supersymmetric
\cite{stanley:on*1, lascoux.schutzenberger:structure,
  fomin.kirillov:yang-baxter, fomin.kirillov:grothendieck}.

\begin{cor} \label{cor:shift}
  Let $w \in S_n$ and $m \leq n$.  Then we have
\[ \begin{split} \Groth_w(c_1,\dots,c_m,a_{m+1},\dots,a_n \,;\,
   c_1,\dots,c_m,b_{m+1},\dots,b_n) = \ \ \ \ \ \ \ \ \ \ \ \ \ \ \ \
   \ \ \ \  \\
   \begin{cases}
   \Groth_u(a_{m+1},\dots,a_n; b_{m+1},\dots,b_n) & \text{if $w =
   1^m\times u$ for some $u$} \\
   0 & \text{otherwise.}
   \end{cases}
\end{split} \]
\end{cor}
\begin{proof}
Apply the theorem to $\Groth_w(c,a,b \,;\, c,b,a)$.
\end{proof}

It was proved in \cite[Prop.~4.1]{feher.rimanyi:schur} that the
Schubert polynomial $\Schub_v(y_v;y)$ is a product of linear factors.
This also follows from our theorem.

\begin{cor} \label{cor:product}
  For $v \in S_n$ we have
\[ \Groth_v(b_v; b) = \prod_{(i,j) \in C(\D_v)} \, 
   \left( 1 - \frac{b_i}{b_{v(j)}} \right) \,.
\]
\end{cor}

Corollaries \ref{cor:bruhat} through \ref{cor:product} above include
all the facts about Schubert polynomials required in \cite{warsawa}.
We remark that Cor.~3 and Prop.~9 of \cite{buch:alternating} are also
special cases of our theorem.

\section{Proof of the main theorem}

Let $R$ be the ring of Laurent polynomials in the variables $a_i$ and
$b_i$, $1 \leq i \leq n$.  For $c \in R$ we set $h_i(c) = (1 +
(1-c)\,s_i) \in \Hecke \otimes R$.  As observed in
\cite{fomin.kirillov:grothendieck}, these elements satisfy the
Yang-Baxter identities $h_i(c)\, h_j(d) = h_j(d)\, h_i(c)$ for $|i-j|
\geq 2$; $h_i(c)\, h_i(d) = h_i(cd)$; and $h_i(c)\, h_{i+1}(cd)\,
h_i(d) = h_{i+1}(d)\, h_i(cd)\, h_{i+1}(c)$.

For $p \geq q$ we furthermore set
\[ A_p^q(c;k) \ = \  h_{k-1+p}(b_p/c) \, h_{k-1+p-1}(b_{p-1}/c) \cdots
   h_{k-1+q}(b_q/c) 
\]
and define, following \cite{fomin.stanley:schubert} and
\cite[(2.1)]{fomin.kirillov:grothendieck}, the product
\[ \Groth^{(n)}(a;b) = A_{n-1}^1(a_1;1) \, A_{n-2}^1(a_2;2) \cdots
   A_1^1(a_{n-1};n-1) \ \in \ \Hecke \otimes R \,.
\]

Fomin and Kirillov have proved that the coefficient of each
permutation $w \in S_n$ in $\Groth^{(n)}(a;b)$ is equal to the
Grothendieck polynomial $\Groth_w(a;b)$ (see Thm.~2.3 and the remark
on page 7 of \cite{fomin.kirillov:grothendieck}, and use the change of
variables $x_i = 1-a_i^{-1}$ and $y_i = 1 - b_i$.)  We claim that the
specialization $\Groth^{(n)}(b_v;b)$ is equal to the south-west to
north-east product of the factors $h_{\nu(i,j)}(b_i/b_{v(j)})$ for all
$(i,j) \in C(\D_v)$.

By descending induction on $q$, the above Yang-Baxter identities imply
that 
\[ A_p^q(c;k-1)\, A_{p-1}^q(d;k)\, h_{q+k-2}(c/d) = A_p^q(d;k-1)\,
   A_{p-1}^q(c;k) \,,
\]
from which we deduce that, for $2 \leq k \leq n-p$, we have
\[ A^{p+1}_{n-k+1}(b_p;k-1) \, A_{n-k}^1(a_k;k) \ = \  
   A_{n-k+1}^{p+1}(a_k;k-1) \, A_{p-1}^1(a_k;k) \,
   A_{n-k}^{p+1}(b_p;k) \,.
\]
By using this identity repeatedly, and setting $\wt a =
(a_2,\dots,a_n)$ and $\wt b = (b_1,\dots,b_{p-1},\\
b_{p+1},\dots,b_n)$, we obtain that
\[ \begin{split}
& \Groth^{(n)}(b_p,\wt a \,;\, b) 
\ \ = \ \ A_{p-1}^1(b_p;1) \, A_{n-1}^{p+1}(b_p;1) \,
  \prod_{k=2}^{n-1} A_{n-k}^1(a_k;k) \\
& \ \ \ = A_{p-1}^1(b_p;1) 
  \left(\prod_{k=2}^{n-p} A_{n-k+1}^{p+1}(a_k;k-1)\, A_{p-1}^1(a_k;k)
  \right) \left( \prod_{k=n-p+1}^{n-1} A_{n-k}^1(a_k;k) \right) \\
& \ \ \ = A_{p-1}^1(b_p;1) 
  \left( 1 \times \Groth^{(n-1)}(\wt a; \wt b) \right) \,.
\end{split} \]
Here ``$1 \times$'' is the operator on $\Hecke \otimes R$ which maps
$s_i$ to $s_{i+1}$ for all $i$.  By setting $p = v(1)$ and $\wt a =
(b_{v(2)},\dots,b_{v(n)})$, the above claim follows by induction, and
our theorem is an immediate consequence of the claim.


\providecommand{\bysame}{\leavevmode\hbox to3em{\hrulefill}\thinspace}
\providecommand{\MR}{\relax\ifhmode\unskip\space\fi MR }
\providecommand{\MRhref}[2]{%
  \href{http://www.ams.org/mathscinet-getitem?mr=#1}{#2}
}
\providecommand{\href}[2]{#2}

\end{document}